\documentclass{elsart}

\usepackage{epsf}

\usepackage{amssymb}
\usepackage{natbib}

\def\SX{{\mathcal X}}

\def\tr{\mbox{\rm trace}}
\def\xb{\mathbf{x}}
\def\Ib{{\mathbf I}}
\def\TT{^\top}
\newtheorem{theo}{\bf Theorem}
\newcommand{\carre} {\hfill \rule{2mm}{2mm} }
\def\me{\epsilon}

\def\SN{{\mathcal N}}

\begin{document}

\begin{frontmatter}


\title{Improvements on removing non-optimal support points in $D$-optimum design algorithms} 


\author[label1]{Radoslav Harman\thanksref{grantRH}}
\thanks[grantRH]{The work of the first author was supported by the VEGA grant No.
1/0264/03 of the Slovak Agency.} \ead{harman@fmph.uniba.sk} and
\author[label2]{Luc Pronzato\corauthref{cor1}\thanksref{grantLP}}
\ead{pronzato@i3s.unice.fr}
\address[label1]{Department of Applied Mathematics and Statistics\\
Faculty of Mathematics, Physics \& Informatics\\
Comenius University, Mlynsk\'a Dolina, 84248 Bratislava, Slovakia}
\address[label2]{Laboratoire I3S, CNRS/UNSA, B\^at.\ Euclide, Les Algorithmes, \\
2000 route des Lucioles, BP 121, 06903 Sophia-Antipolis Cedex, France}
\corauth[cor1]{author for correspondence}
\thanks[grantLP]{The work of the second author was
partially supported by the IST Programme of the European
Community, under the PASCAL Network of Excellence,
IST-2002-506778. This publication only reflects the authors view.}


\begin{abstract}
We improve the inequality used in \citep{Pa03} to remove points
from the design space during the search for a $D$-optimum design.
Let $\xi$ be any design on a compact space $\mathcal{X} \subset
\mathbb{R}^m$ with a nonsingular information matrix, and let
$m+\epsilon$ be the maximum of the variance function
$d(\xi,\mathbf{x})$ over all $\mathbf{x} \in \mathcal{X}$. We
prove that any support point $\mathbf{x}_{*}$ of a $D$-optimum
design on $\mathcal{X}$ must satisfy the inequality
$d(\xi,\mathbf{x}_{*}) \geq
m(1+\epsilon/2-\sqrt{\epsilon(4+\epsilon-4/m)}/2)$. We show that
this new lower bound on $d(\xi,\mathbf{x}_{*})$ is, in a sense,
the best possible, and how it can be used to accelerate algorithms
for $D$-optimum design.
\end{abstract}

\end{frontmatter}

\section{Introduction}
Let $\mathcal{X} \subseteq \mathbb{R}^m$ be a compact design space
and let $\Xi$ be the set of all designs (i.e., finitely supported
probability measures) on $\mathcal{X}$. For any $\xi \in \Xi$, let
$$
 \mathbf{M}(\xi)=\int_{\mathcal{X}}\mathbf{x}\mathbf{x}\TT \; \xi(\mathbf{dx})
$$
denote the information matrix. Suppose that there exists a design
with nonsingular information matrix and let $\Xi^{+}$ be the set
of such designs. Let $\xi^{*}$ denote a $D$-optimum design, that
is, a measure in $\Xi$ that maximizes $\det\mathbf{M}(\xi)$, see,
e.g., \citep{Fedorov72}. Note that a $D$-optimum design always
exists and that the $D$-optimum information matrix
$\mathbf{M}_{*}=\mathbf{M}(\xi^{*})$ is unique. For any $\xi \in
\Xi^{+}$ denote $d(\xi,\cdot):\mathcal{X} \to [0,\infty)$ the
variance function defined by
$$
d(\xi,\mathbf{x})=\mathbf{x}\TT \mathbf{M}^{-1}(\xi)\mathbf{x} \,.
$$

The celebrated Kiefer-Wolfowitz Equivalence Theorem (1960)\nocite{KieferW60} writes as follows.

\begin{theo}\label{Thm:EQ}
The following three statements are equivalent:
\begin{description}
    \item[(i)] $\xi^*$ is $D$-optimum;
    \item[(ii)] $\max_{\xb\in\SX} d(\xi^*,\xb) =m$;
    \item[(iii)] $\xi^*$ minimizes $\max_{\xb\in\SX} d(\xi,\xb)$,
    $\xi\in\Xi^+$.
\end{description}
\end{theo}
Notice that
$$
\int_{\mathcal{X}} d(\xi^*,\mathbf{x}) \; \xi^{*}(d\mathbf{x}) =
    \int_{\mathcal{X}}\mathbf{x}\TT \mathbf{M}_{*}^{-1}\mathbf{x} \; \xi^{*}(d\mathbf{x})=
    \tr(\mathbf{M}_{*}\mathbf{M}_{*}^{-1})=m\,.
$$
Hence, (ii) of Theorem \ref{Thm:EQ} implies that for any support
point $\mathbf{x}_{*}$ of the design $\xi^{*}$ (i.e., for a point
satisfying $\xi^{*}(\mathbf{x}_{*})>0$), we have
 \begin{equation} \label{Eq:SUP}
   d(\xi^*,\mathbf{x}_*)=m \,.
 \end{equation}
 In the next section we show that the equality (\ref{Eq:SUP}) can be used to prove that
 $$
 \forall \xi \in \Xi^{+}\,, \ d(\xi,\mathbf{x}_{*}) \geq m\lambda_1^*(\xi)
 $$
 where $\lambda_1^*$ depends on $\xi$ only via the maximum of $d(\xi,\cdot)$ over the design space $\SX$.
 Hence, we can test candidate support points by using any finite number of design measures
 $\xi \in \Xi^{+}$, e.g., those that are generated by a design algorithm on its way towards the
 optimum: any point that does not pass the test defined by $\xi^k$ of iteration $k$ need not be
 considered for further investigations and can thus be removed from the design
 space.
\section{A necessary condition for candidate support points}

For $\xi$ a design in $\Xi^{+}$ denote
$\mathbf{M}=\mathbf{M}(\xi)$,
$$
\mathbf{H}=\mathbf{M}^{-1/2}\mathbf{M}_{*}\mathbf{M}^{-1/2}
$$
and $\lambda_1 \leq \lambda_2 \leq \cdots \leq \lambda_m$ the
eigenvalues of $\mathbf{H}$. Notice that $\lambda_1>0$ and that
the eigenvalues depend on the design $\xi$ as well as on the
$D$-optimum information matrix $\mathbf{M}_{*}$. Let
$\mathbf{x}_{*}$ be a support point of a $D$-optimum design and
let
$\mathbf{y}_{*}=\mathbf{H}^{-1/2}\mathbf{M}^{-1/2}\mathbf{x}_{*}$.
The equality (\ref{Eq:SUP}) can be written in the form
$\mathbf{y}_{*}\TT \mathbf{y}_{*} = m$ which implies:
\begin{equation} \label{Ineq:MEGN}
   d(\xi,\mathbf{x}_{*})=\mathbf{x}_{*}\TT \mathbf{M}^{-1}\mathbf{x}_{*}=
   \mathbf{y}_{*}\TT \mathbf{H}\mathbf{y}_{*} \geq
   \lambda_1 \mathbf{y}_{*}\TT \mathbf{y}_{*} = m \lambda_1 \,.
\end{equation}
To be able to use the inequality (\ref{Ineq:MEGN}), we need to
derive a lower bound $\lambda_1^*$ on $\lambda_1$ that does not
depend on the unknown matrix $\mathbf{M}_{*}$.

Theorem \ref{Thm:EQ}-(ii) implies
\begin{eqnarray*}
   \sum_{i=1}^{m}\lambda_i^{-1} &=& \tr(\mathbf{H}^{-1})\\
   &=& \tr(\mathbf{M}_{*}^{-1}\mathbf{M})=
   \int_{\mathcal{X}}\mathbf{x}\TT \mathbf{M}_{*}^{-1}\mathbf{x} \; \xi(d\mathbf{x}) =
   \int_{\mathcal{X}}d(\xi^*,\mathbf{x}) \; \xi(d\mathbf{x}) \leq
   m\,.
\end{eqnarray*}
Also,
\begin{eqnarray*}
   \sum_{i=1}^{m}\lambda_i &=& \tr(\mathbf{H}) \\
   &=& \tr(\mathbf{M}_{*}\mathbf{M}^{-1}) =
   \int_{\mathcal{X}}\mathbf{x}\TT \mathbf{M}^{-1}\mathbf{x} \; \xi^*(d\mathbf{x}) \leq
   \max_{\mathbf{x} \in \mathcal{X}}\mathbf{x}\TT \mathbf{M}^{-1}\mathbf{x} =
   m+\epsilon\,,
\end{eqnarray*}
where we used the notation
\begin{equation}\label{epsilon}
    \epsilon=\epsilon(\xi)=\max_{\mathbf{x} \in \mathcal{X}}\mathbf{x}\TT
\mathbf{M}^{-1}\mathbf{x}-m \geq 0 \,.
\end{equation}
For $m=1$ we directly obtain the lower bound $\lambda_1 \geq
\lambda_1^{*}=1$. For $m>1$, the Lagrangian for the minimisation
of $\lambda_1$ subject to $\sum_{i=1}^{m}\lambda_i^{-1} \leq m$
and $\sum_{i=1}^{m}\lambda_i \leq m+\epsilon$ is given by
$$
 \mathcal{L}(\lambda,\mu_1,\mu_2)=\lambda_1+
 \mu_1\left(\sum_{i=1}^{m}\lambda_i^{-1} - m\right) +
 \mu_2\left(\sum_{i=1}^{m}\lambda_i - m -\epsilon \right)
$$
with $\lambda=(\lambda_1,...,\lambda_m)\TT$, $\mu_1,\mu_2 \geq 0$.
The stationarity of $\mathcal{L}(\lambda,\mu_1,\mu_2)$ with
respect to the $\lambda_i$'s and the Kuhn-Tucker conditions
$$
 \mu_1\left(\sum_{i=1}^{m}\lambda_i^{-1} - m\right)=0 \mbox{, }
 \mu_2\left(\sum_{i=1}^{m}\lambda_i - m -\epsilon \right) = 0
$$
give $\lambda_i=L$ for $i=2,...,m$, with $\lambda_1$ and $L$
satisfying
\begin{eqnarray*}
  && \lambda_1^{-1}+(m-1)L^{-1}=m \\
  && \lambda_1+(m-1)L=m+\epsilon \,.
\end{eqnarray*}
The solution is thus
\begin{equation} \label{Ineq:LB}
  \lambda_1^{*}=1+\frac{\epsilon}{2}-\frac{\sqrt{\epsilon(4+\epsilon-4/m)}}{2} \leq 1
\end{equation}
and $\lambda_i^{*}=L^{*}=(m-1)/(m-1/\lambda_1^{*}) \geq 1$,
$i=2,...,m$. Notice that the bound (\ref{Ineq:LB}) gives
$\lambda_1^{*}=1$ when $m=1$ and can thus be used for any
dimension $m\geq 1$. By substituting $\lambda_1^*$ for $\lambda_1$
in (\ref{Ineq:MEGN}) we obtain the following result.

\begin{theo}
 For any design $\xi \in \Xi^+$, any point $\mathbf{x}_* \in \mathcal{X}$ such that
\begin{equation}\label{**}
     d(\xi,\mathbf{x}_*) <
 h_m(\epsilon)= m \left[ 1+\frac{\epsilon}{2}-\frac{\sqrt{\epsilon(4+\epsilon-4/m)}}{2} \right]
\end{equation}
where $\epsilon=\max_{\mathbf{x} \in \mathcal{X}}
d(\xi,\mathbf{x}) -m$, cannot be a support point of a D-optimum
design measure.
\end{theo}

The inequality in \citep{Pa03} uses
\begin{equation}\label{tildeh}
   \tilde{h}_m(\epsilon)=m \left[ 1+\frac{\epsilon}{2}-\frac{\sqrt{\epsilon(4+\epsilon)}}{2}
   \right]\,.
\end{equation}
Notice, that $m\geq h_m(\epsilon)>\tilde{h}_m(\epsilon)$ for all
integer $m\geq 1$ and all $\epsilon>0$, and that $\lim_{\epsilon
\to \infty}h_m(\epsilon)=1$ while $\lim_{\epsilon \to
\infty}\tilde{h}_m(\epsilon)=0$. The new bound is thus always
stronger, especially for large values of $\epsilon$, i.e. when the
design $\xi$ is far from being optimum. Although in practice the
improvement over (\ref{tildeh}) can be marginal, see the example
below, the important result here is that the bound (\ref{**})
cannot be improved. Indeed, when $m=1$, $h_1(\me)=1$ for any
$\me>0$ which is clearly the best possible bound. When $m\geq 2$,
$h_m(\epsilon)$ is the tightest lower bound on the variance
function $d(\xi,\mathbf{x}_{*})$ at a $D$-optimal support point
$\xb_*$ that depends only on $m$ and $\epsilon$, in the sense of
the following theorem.

\begin{theo}
For any integer $m\geq 2$ and any $\epsilon,\delta>0$ there exist
a compact design space $\mathcal{X} \subset \mathbb{R}^m$, a
design $\xi$ on $\mathcal{X}$ and a point $\mathbf{x}_{*} \in
\mathcal{X}$ supporting a $D$-optimum design on $\SX$ such that
$\epsilon=\max_{\xb \in \mathcal{X}} d(\xi,\mathbf{x})-m$ and
$$
d(\xi,\mathbf{x}_{*})<h_m(\epsilon)+\delta \,.
$$
\end{theo}
{\em Proof.} Denote $h=h_m(\epsilon)$ and $k=2^{m-1}$. Let
$\mathbf{x}_1,...,\mathbf{x}_k$ correspond to the $k$ vectors of
$\mathbb{R}^m$ of the form
$$
\left(\sqrt{\frac{1}{h}}, \pm \sqrt{\frac{h-1}{h(m-1)}},\ldots,\pm
\sqrt{\frac{h-1}{h(m-1)}} \right)\TT
$$
and let $\mathbf{y}_1,\ldots,\mathbf{y}_k$ correspond to the $k$
vectors $\left(\sqrt{1/m}, \pm \sqrt{1/m},\ldots,\pm
\sqrt{1/m}\right)\TT$. Take
$\mathbf{x}_{*}=(\sqrt{b},0,\ldots,0)\TT \in \mathbb{R}^m$ with
$1<b<\min\left\{(\epsilon+m)/h, (h+\delta)/h\right\}$,
$\mathcal{X}$ as the finite set
$\mathcal{X}=\{\mathbf{x}_{1},...,\mathbf{x}_{k},
  \mathbf{y}_{1},...,\mathbf{y}_{k},\mathbf{x}_{*}\}$
and let $\xi$ be the uniform probability measure on
$\mathbf{x}_{1},...,\mathbf{x}_{k}$. Note that $\mathbf{M}(\xi)$
is a diagonal matrix with diagonal elements $\left(1/h,
(h-1)/[h(m-1)],\ldots,(h-1)/[h(m-1)]\right)$. One can easily
verify that
$$
\max_{\mathbf{x} \in \mathcal{X}}\mathbf{x}\TT
\mathbf{M}^{-1}(\xi)\mathbf{x}-m = \me \mbox{ and }
d(\xi,\mathbf{x}_{*})=\mathbf{x}_{*}\TT
\mathbf{M}^{-1}(\xi)\mathbf{x}_{*} =b\,h<h_m(\epsilon)+\delta \,.
$$
The uniform probability measure $\eta$ on
$\mathbf{y}_1,...,\mathbf{y}_k$ is $D$-optimum on
$\mathcal{X}/\{\mathbf{x}_{*}\}$, as can be directly verified by
checking (ii) of the Equivalence Theorem \ref{Thm:EQ}. On the
other hand, $\eta$ is not $D$-optimum on $\mathcal{X}$ since
$\mathbf{x}_{*}\TT \mathbf{M}^{-1}(\eta)\mathbf{x}_{*}=b\,m>m$,
which implies that $\mathbf{x}_{*}$ must support a $D$-optimum
design on $\mathcal{X}$. \carre

\vspace{0.3cm} {\bf Example:} We consider a series of problems
defined by the construction of the minimum covering ellipse for an
initial set of 1000 random points in the plane, i.i.d.\
$\SN(0,\Ib_2)$. These problems correspond to $D$-optimum design
problems for randomly generated $\SX \subset \mathbb{R}^3$, see
\citet{Titterington75,Titterington78}. The following recursion can
thus be used:\begin{equation}\label{recursion}
    w_i^{k+1} = w_i^k \frac{d(\xi^k,\xb_i)}{m} \,, \
    i=1,\ldots,q(k)\,,
\end{equation}
where $k\geq 0$, $w_i^k=\xi^k(\xb_i)$ is the weight given by the
discrete design $\xi^k$ to the point $\xb_i$ and $q(k)$ is the
cardinality of $\SX$ at iteration $k$. In the original algorithm,
$q(k)=q(0)$ for all $k$ and, initialized at a $\xi^0$ that gives a
positive weight at each point of $\SX$, the algorithm converges
monotonically to the optimum, see \citep{Torsney83} and
\citep{Titterington76}. The tests (\ref{**}) and (\ref{tildeh})
can be used to decrease $q(k)$: at iteration $k$, any design point
$\xb_j$ satisfying $d(\xi^k,\xb_j)< h_m[\me(\xi^k)]$, see
(\ref{epsilon}, \ref{**}), or $d(\xi^k,\xb_j)< \tilde
h_m[\me(\xi^k)]$, see (\ref{epsilon}, \ref{tildeh}), can be
removed from $\SX$. The total weight of the points that are
cancelled is then reallocated to the $\xb_i$'s that stay in $\SX$
(e.g., proportionally to $w_i^k$).

\begin{figure}[ht]
\centerline{\vbox{\epsfxsize=8cm\epsfbox{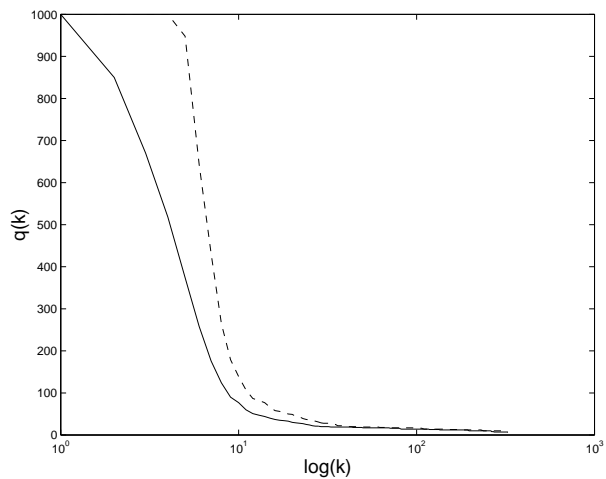}}}
\caption{$q(k)$ as a function of $\log(k)$: cancellation based on
(\ref{**}) in solid line, on (\ref{tildeh}) in dashed line.}
\label{F:ellips}
\end{figure}

Figure \ref{F:ellips} presents a typical evolution of $q(k)$ as a
function of $\log(k)$ for $\xi^0$ uniform on $\SX$ and shows the
superiority of the test (\ref{**}) over (\ref{tildeh}). The
improvement is especially important in the first iterations, when
the design $\xi^k$ is far from the optimum. Define $k^*(\delta)$
as the number of iterations required to reach a given precision
$\delta$,
$$
k^*(\delta) = \min\left\{k\geq 0: \me(\xi^k) < \delta \right\} \,,
$$
with $\me(\xi^k)$ defined by (\ref{epsilon}). Notice that from the
concavity of $\log\det \mathbf{M}(\xi)$ we have
\begin{eqnarray*}
  \log \det \mathbf{M}(\xi^*) - \log\det
\mathbf{M}(\xi^{k^*(\delta)}) &\leq&  \frac{\partial\log \det
\mathbf{M}[(1-\alpha)\xi^{k^*(\delta)}+\alpha
\xi^*]}{\partial\alpha}_{|\alpha=0} \\
  &&= \int_{\SX} d(\xi^{k^*(\delta)},\xb)\, \xi^*(d\xb) - m < \delta \,.
\end{eqnarray*}
Table \ref{TB:ellips} shows the influence on the algorithm
(\ref{recursion}) of the cancellation of points based on the tests
(\ref{**}) and (\ref{tildeh}), in terms of $k^*(\delta)$, of the
corresponding computing time $T(\delta)$, the number of support
points $n(\delta)$ of $\xi^{k^*(\delta)}$ and the first iteration
$k_{10}$ when $\xi^k$ has 10 support points or less, with
$\delta=10^{-3}$. The results are averaged over 1000 independent
problems. The values of $k^*(\delta)$ and $k_{10}$ are rounded to
the nearest larger integer, the computing time for the algorithm
with the cancellation of points based on (\ref{**}) is taken as
reference and set to 1 (the algorithm without cancellation was at
least 4.5 times slower in all the 1000 repetitions). Although
cancelling points has little influence on the number of iterations
$k^*(\delta)$, is renders the iterations simpler: on average the
introduction of the test (\ref{**}) in the algorithm
(\ref{recursion}) makes it about 30 times faster.

\begin{table}[ht]
\centering\begin{tabular}{c|cccc} \hline
Algorithm & $k^*(\delta)$ & $T(\delta)$ & $n(\delta)$ & $k_{10}$ \\
(\ref{recursion}) & 252 & 31.6 & 1000 & $-$ \\
(\ref{recursion}) and (\ref{tildeh}) & 248 & 1.4 & 5.8 & 82 \\
(\ref{recursion}) and (\ref{**}) & 247 & 1 & 5.5 & 66 \\
 \hline
\end{tabular}
\caption{Influence of the tests (\ref{**}) and (\ref{tildeh}) on
the average performance of the algorithm (\ref{recursion}) for the
minimum covering ellipse problem (1000 repetitions,
$\delta=10^{-3}$).} \label{TB:ellips}
\end{table}

The influence of the cancellation on the performance of the
algorithm can be further improved as follows. Let $(k_j)_j$ denote
the subsequence corresponding to the iterations where some points
are removed from $\SX$. We have $j \leq q(0)$, the cardinality of
the initial $\SX$, and the convergence of the algorithm
(\ref{recursion}) is therefore maintained whatever the heuristic
rule used at the iterations $k_j$ for updating the weights of the
points that stay in $\SX$ (provided these weights remain strictly
positive). The following one has been found particularly efficient
on a series of examples: for all $t\in T_j$, the set of indices
corresponding to the points that stay in $\SX$ at iteration $k_j$,
replace $w_t^{k_j}$ by
$$
    {w'_t}^{k_j} =  \frac{z_t}{\sum_{s \in T_j}z_s}
\ \mbox{ where} \ z_t = \left\{ \begin{array}{ll} A w_t^{k_j} & \mbox{ if } d(\xi^{k_j},\xb_t) \geq m\\
w_t^{k_j} & \mbox{ otherwise } \end{array} \right.
$$
for some $A \geq 1$. A final remark is that by including the test
(\ref{**}) in the algorithm (\ref{recursion}) one can in general
quickly identify potential support points for an optimum design.
When the number $n$ of these points is small enough, switching to
a more standard convex-programming algorithm for the optimization
of the $n$ associated weights might then form a very efficient
strategy.



\end{document}